# Recursive Harmonic Numbers and Binomial Coefficients

Aung Phone Maw and Aung Kyaw
Department of Mathematics
University of Yangon
Yangon, 11041, Myanmar
uaungkyaw70@gmail.com

**Abstract.** We define recursive harmonic numbers as a generalization of harmonic numbers. The table of recursive harmonic numbers, which is like Pascal's triangle, is constructed. A formula for recursive harmonic numbers containing binomial coefficients is also presented.

## 1. Recursive Harmonic Numbers

For a positive integer $n$, a harmonic number $H_n$ is defined as $H_n = \sum_{k=1}^{n} \frac{1}{k}$. Here we define $m^{\text{th}}$ recursive harmonic number $H_n^{(m)}$ as follows:

$$H_n^{(0)} = 1 \; ; \; H_n^{(m)} = \sum_{k=1}^{n} \frac{1}{k} H_k^{(m-1)} \text{ for any positive integer } m.$$

Since $H_n = \sum_{k=1}^{n} \frac{1}{k} = \sum_{k=1}^{n} \frac{1}{k} \cdot 1 = \sum_{k=1}^{n} \frac{1}{k} H_k^{(0)} = H_n^{(1)}$, one can see that $m^{\text{th}}$ recursive harmonic number is a generalization of a harmonic number.

## 2. Table of Recursive Harmonic Numbers

From the definition of $m^{\text{th}}$ recursive harmonic number, $H_n^{(0)} = 1$ and $H_1^{(m)} = 1$. For every $n \geq 2$ we have

$$H_n^{(m)} = \sum_{k=1}^{n} \frac{1}{k} H_k^{(m-1)}$$

$$= \sum_{k=1}^{n-1} \frac{1}{k} H_k^{(m-1)} + \frac{1}{n} H_n^{(m-1)}$$

$$H_n^{(m)} = H_{n-1}^{(m)} + \frac{1}{n} H_n^{(m-1)}$$

From these facts we can construct the table of recursive harmonic numbers like Pascal's triangle as follows:

| $n$ \ $m$ | 0 | 1 | 2 | 3 | 4 |
|---|---|---|---|---|---|
| 1 | 1 | 1 | 1 | 1 | 1 |
| 2 | 1 | $\frac{3}{2}$ | $\frac{7}{4}$ | $\frac{15}{8}$ | $\frac{31}{16}$ |
| 3 | 1 | $\frac{11}{6}$ | $\frac{85}{36}$ | $\frac{575}{216}$ | $\frac{3661}{1296}$ |
| 4 | 1 | $\frac{25}{12}$ | $\frac{415}{144}$ | $\frac{5845}{1728}$ | $\frac{76111}{20736}$ |

### 3. Recursive Harmonic Numbers and Binomial Coefficients

A well-known formula for harmonic numbers containing binomial coefficients is

$$H_n = H_n^{(1)} = \sum_{k=1}^{n} (-1)^{k+1} \frac{1}{k} \binom{n}{k}.$$

We will show that

$$H_n^{(m)} = \sum_{k=1}^{n} (-1)^{k+1} \frac{1}{k^m} \binom{n}{k}.$$

**Proof.** We will prove by induction.

When $n = 1$, $H_1^{(m)} = 1 = \sum_{k=1}^{1} (-1)^{k+1} \frac{1}{k^m} \binom{1}{k}$.

When $m = 0$, $H_n^{(0)} = 1$ and

$$\sum_{k=1}^{n} (-1)^{k+1} \frac{1}{k^0} \binom{n}{k} = \binom{n}{1} - \binom{n}{2} + \binom{n}{3} - \binom{n}{4} + \cdots + (-1)^{n+1} \binom{n}{n} = 1.$$

Therefore $H_n^{(0)} = \sum_{k=1}^{n} (-1)^{k+1} \frac{1}{k^0} \binom{n}{k}$.

By assuming the formula is true for $H_n^{(m-1)}$ and $H_{n-1}^{(m)}$, we will show that the formula is true for $H_n^{(m)}$, $n \geq 2, m \geq 1$. Since $H_n^{(m)} = H_{n-1}^{(m)} + \frac{1}{n} H_n^{(m-1)}$, then

$$H_n^{(m)} = H_{n-1}^{(m)} + \frac{1}{n} H_n^{(m-1)}$$

$$= \sum_{k=1}^{n-1} (-1)^{k+1} \frac{1}{k^m} \binom{n-1}{k} + \frac{1}{n} \sum_{k=1}^{n} (-1)^{k+1} \frac{1}{k^{m-1}} \binom{n}{k}$$

$$= \sum_{k=1}^{n-1} (-1)^{k+1} \frac{1}{k^m} \frac{n-k}{n} \binom{n}{k} + (-1)^{n+1} \frac{1}{n \cdot n^{m-1}} \cdot 1 + \sum_{k=1}^{n-1} (-1)^{k+1} \frac{1}{n \cdot k^{m-1}} \binom{n}{k} +$$

$$= \sum_{k=1}^{n-1} (-1)^{k+1} \frac{1}{k^m} \frac{(n-k)+k}{n} \binom{n}{k} + (-1)^{n+1} \frac{1}{n^m}$$

$$= \sum_{k=1}^{n-1} (-1)^{k+1} \frac{1}{k^m} \binom{n}{k} + (-1)^{n+1} \frac{1}{n^m}$$

$$H_n^{(m)} = \sum_{k=1}^{n} (-1)^{k+1} \frac{1}{k^m} \binom{n}{k}$$

Other formulas involving harmonic numbers and binomial coefficients can be found in [1, 2, 3] and others.

### References

1. J. Choi, Finite Summation Formulas involving Binomial Coefficients, Harmonic Numbers and Generalized Harmonic Numbers, *J. Inequal. Appl.* **2013**:49

2.  W. Chu and Q. Yan, Combinatorial Identities on Binomial Coefficients and Harmonic Numbers, *Util. Math.* 75 (2008) 51-66

3.  M.J. Kronenburg, Some Combinatorial Identities some of which involving Harmonic Numbers, arXiv: 1103.1268v3 [math.CO] 12 Jan 2017